\documentclass[10pt]{amsart}            %%%  {article}

\usepackage{latexsym,amssymb}
\usepackage{amsmath,amsthm}      %%% for better alignments and theorems
\usepackage{amsopn}              %%% for `operator names' like `\det'
\usepackage{cite}                %%% for better handling of citations

\title{Differential forms canonically associated
to even-dimensional compact conformal manifolds}

\author{William J. Ugalde}

\newcommand{\nnn}[1]{\eqref{#1}}

\newtheorem{theorem}{Theorem}[section]
\newtheorem{lemma}[theorem]{Lemma}
\newtheorem{proposition}[theorem]{Proposition}

\newtheorem{definition}[theorem]{Definition}

\newcommand{\bbC}{\mathbb{C}}

\newcommand{\bbR}{\mathbb{R}}

\newcommand{\cH}{\mathcal{H}}

\newcommand{\cP}{\mathcal{P}}

\renewcommand{\a}{\alpha}
\renewcommand{\b}{\beta}
\renewcommand{\d}{\delta}
\newcommand{\ga}{\gamma}

\newcommand{\La}{\Lambda}
\newcommand{\Na}{\nabla}
\newcommand{\s}{\sigma}

\def\<#1,#2>{\langle{#1},{#2}\rangle} %% my version for inner product
\DeclareMathOperator{\tr}{trace}  %%% uses `amsopn' package

\def\set#1{\{\,#1\,\}}          %% set notation
\def\del{\partial}              %% partial derivatives symbol
\def\sideremark#1{\ifvmode\leavevmode\fi\vadjust{\vbox to0pt{\vss
\hbox to 0pt{\hskip\hsize\hskip1em
\vbox{\hsize2cm\tiny\raggedright\pretolerance10000
\noindent #1\hfill}\hss}\vbox to8pt{\vfil}\vss}}}

\def\lcovd#1#2{{#1}_{\,|\,#2}\,}        %% for f|_i
\def\ucovd#1#2{{#1}_{\,|\,}{}^{#2}\,}   %% for f|^i

\def\llcovd#1#2#3{{#1}_{\,|\,#2\,\,#3}\,}             %% for  f|_{ij}
\def\lucovd#1#2#3{{#1}_{\,|\,#2\,}{}^{#3}\,}          %% for  f|_i^j
\def\ulcovd#1#2#3{{#1}_{\,|\,}{}^{#2}{}_{\,#3}\,}     %% for  f|^i_j
\def\uucovd#1#2#3{{#1}_{\,|\,}{}^{#2\,#3}\,}          %% for  f|^{ij}

\def\ulucovd#1#2#3#4{{#1}_{\,|\,}{}^{#2\,}{}_{#3\,}{}^{#4}\,} %% for  f|^i_j^k
\def\llucovd#1#2#3#4{{#1}_{\,|\,#2\,#3\,}{}^{#4}\,}           %% for  f|_{ij}^k
\def\uuucovd#1#2#3#4{{#1}_{\,|\,}{}^{#2\,#3\,#4}\,}           %% for  f|^{ijk}
\def\ullcovd#1#2#3#4{{#1}_{\,|\,}{}^{#2\,}{}_{#3\,#4}\,}      %% for  f|^i_{jk}
\def\lllcovd#1#2#3#4{{#1}_{\,|\,#2\,#3\,#4}\,}                %% for  f|_{ijk}

\def\lllucovd#1#2#3#4#5{{#1}_{\,|\,#2\,#3\,#4\,}{}^{#5}\,}           %% for  f|_{ijk}^l
\def\lulucovd#1#2#3#4#5{#1_{\,|\,#2\,}{}^{#3\,}{}_{#4\,}{}^{#5}\,}   %% for  f|_i^j_k^l
\def\uulucovd#1#2#3#4#5{#1_{\,|\,}{}^{#2\,#3\,}{}_{#4\,}{}^{#5}\,}   %% for  f|^{ij}_k^l

\def\llulucovd#1#2#3#4#5#6{#1_{\,|\,#2\,#3\,}{}^{#4\,}{}_{#5\,}{}^{#6}\,}
\def\ululucovd#1#2#3#4#5#6{{#1}_{\,|\,}{}^{#2\,}{}_{#3\,}{}^{#4\,}{}_{#5\,}{}^{#6}\,}
\def\lullucovd#1#2#3#4#5#6{{#1}_{\,|\,#2\,}{}^{#3\,}{}_{#4\,#5}{}^{#6}\,}
\renewcommand{\lcovd}[2]{{#1}_{\,;\,#2}}     %% f|_i
\renewcommand{\ucovd}[2]{{#1}_{\,;}{}^{#2}}  %% f|^i

\renewcommand{\llcovd}[3]{{#1}_{\,;\,#2#3}}         %% f|_{ij}
\renewcommand{\lucovd}[3]{{#1}_{\,;\,#2}{}^{#3}}    %% f|_i^j
\renewcommand{\ulcovd}[3]{{#1}_{\,;}{}^{#2}{}_{#3}} %% f|^i_j
\renewcommand{\uucovd}[3]{{#1}_{\,;}{}^{#2#3}}      %% f|^{ij}

\renewcommand{\ulucovd}[4]{{#1}_{\,;}{}^{#2}{}_{#3}{}^{#4}} %% f|^i_j^k
\renewcommand{\llucovd}[4]{{#1}_{\,;\,#2#3}{}^{#4}}         %% f|_{ij}^k
\renewcommand{\uuucovd}[4]{{#1}_{\,;}{}^{#2#3#4}}           %% f|^{ijk}
\renewcommand{\ullcovd}[4]{{#1}_{\,;}{}^{#2}{}_{#3#4}}      %% f|^i_{jk}
\renewcommand{\lllcovd}[4]{{#1}_{\,;\,#2#3#4}}              %% f|_{ijk}

\renewcommand{\lllucovd}[5]{{#1}_{\,;\,#2#3#4}{}^{#5}}         %% f|_{ijk}^l
\renewcommand{\lulucovd}[5]{#1_{\,;\,#2}{}^{#3}{}_{#4}{}^{#5}} %% f|_i^j_k^l
\renewcommand{\uulucovd}[5]{#1_{\,;}{}^{#2#3}{}_{#4}{}^{#5}}   %% f|^{ij}_k^l 

\renewcommand{\llulucovd}[6]{#1_{\,;\,#2#3}{}^{#4}{}_{#5}{}^{#6}}
\renewcommand{\ululucovd}[6]{{#1}_{\,;}{}^{#2}{}_{#3}{}^{#4}{}_{#5}{}^{#6}}
\renewcommand{\lullucovd}[6] {{#1}_{\,;\,#2}{}^{#3}{}_{#4#5}{}^{#6}}
\newcommand{\eps}{\varepsilon}      %% abbreviation
\renewcommand{\flat}{\mathrm{flat}} %% subscript for flat case
\newcommand{\confflat}{\mathrm{conf.\,flat}} %% conformally flat case   

\DeclareMathOperator{\dvol}{dvol} %% Riemannian volume form
\DeclareMathOperator{\Hess}{Hess} %% Hessian
\DeclareMathOperator{\Rc}{Rc}     %% Ricci tensor
\DeclareMathOperator{\Sc}{Sc}     %% scalar curvature
\DeclareMathOperator{\Wres}{Wres} %% Wodzicki residue
\DeclareMathOperator{\wres}{wres} %% Wodzicki residue density

\begin{document}

\begin{abstract}  
On a 6-dimensional, conformal, oriented, compact  manifold $M$ without boundary, 
we compute a whole family of differential forms $\Omega_6(f,h)$ of order 6, with $f,h \in C^\infty(M).$
Each of these 
forms will be symmetric on $f,$ and $h,$ conformally invariant, and such that $\int_M f_0 \Omega_6(f_1,f_2)$ 
defines a Hochschild 2-cocycle over the algebra $C^\infty(M).$  In the particular 
6-dimensional conformally flat case, we compute the unique one satisfying 
$\Wres(f_0[F,f][F,h]) = \int_M f_0\Omega_6(f,h)$ for $(\cH,F)$ the Fredholm module associated by 
A. Connes \cite{Con1} to the manifold $M,$ and $\Wres$ the Wodzicki residue.\par  
\end{abstract}

%% 53A30 = Conformal differential geometry
%% 35S05 = General theory of PsDO
%% 58J40 = Pseudodifferential and Fourier integral operators on manifolds
%% 58B34 = Noncommutative geometry (à la Connes)
%% 46L87 = Noncommutative differential geometry

\maketitle  %%% Outputs the title, author(s) and date

\tableofcontents

\pagestyle{myheadings}
\markboth{William \ J. \ Ugalde}{Differential forms associated to
conformal manifolds}

{ {\renewcommand{\thefootnote}{\fnsymbol{footnote}}}}
\footnotetext{Expanded version of a talk given at ``The 6th. Conference on Cliffod algebras and their 
applications in Mathematical Physics.''   Tennessee Technological University,   Cookeville,  
Tennessee,  May 20-25, 2002.}

%%%%%%%%%%%%%%%%start body of the paper here

%%%%%%%%%%%%%%%%%%%%%%%%  INTRODUCTION %%%%%%%%%%%%%%%%%%%%%%%%%%%%%%%%%%%%%%%

\section{Introduction}

For a compact, oriented manifold $M$ of even dimension $n=2l,$ endowed with a conformal 
structure,
there is a canonically associated Fredholm module
$(\cH,F)$ \cite{Con1}.  $\cH$ is the 
Hilbert space of square integrable forms of middle dimension 
$$\cH=L^2(M,\Lambda^{n/2}_{\bbC} T^*M),$$
in which functions on $M$ act as
multiplication operators. $F$ is the pseudodifferential operator of order 0 acting in $\cH,$
obtained from the orthogonal projection $P$ on the image of $d,$ by the relation $F=2P-1.$  
From the Hodge decomposition theorem \cite{Warner} 
it is easy to see that $F$ preserves the finite 
dimensional space of harmonic forms $H^l$, and that $F$ restricted to the $\cH\ominus H^l$
is given by
\begin{equation}
\label{defofF}
F=\frac{d \d-\d d}{d \d+\d d},
\end{equation}
in terms of a Riemannian metric compatible with the conformal structure of $M.$
Both $\cH$ and $F$ are independent of the metric in the conformal class \cite[Section~IV.4.$\ga$]{Koran}.

By considering a Riemann surface $\Sigma,$ and by considering instead of $d X$ its 
quantized version $[F,X],$ Connes (Chapter~4~\cite{Koran}) quantized the Polyakov 
action
as the Dixmier trace of the operator $\eta_{ij} d X^i d X^j$
\begin{equation}
\label{quantPolyakov}
\frac{1}{2 \pi} \int_\Sigma \eta_{ij} d X^i \wedge \star d X^j = 
-\frac{1}{2} \tr_\omega\bigl(\eta_{ij}[F,X^i][F,X^j]\bigr).
\end{equation}
Because of the fact that the Wodzicki residue extends uniquely the Dixmier trace as a
trace on the algebra of pseudodifferential operators \cite{Wo}, this quantized Polyakov 
action has sense in the general even-dimensional case.
Because of the Connes' trace theorem (see Theorem~7.18~~\cite{Polaris}), 
the Dixmier trace and the Wodzicki residue of an elliptic pseudodifferential operator 
of order $-n$ in 
an $n$-dimensional manifold are proportional by a factor of $n(2\pi)^n.$  
In the 2-dimensional 
case the factor is $8\pi^2$ and so the quantized Polyakov action 
\nnn{quantPolyakov} can be written as, 
\begin{equation}
\label{generalPolyakov}
-16\pi^2 I=\Wres \bigl(\eta_{ij} [F,X^i][F,X^j]\bigr)
\end{equation}
which determines, by using the general formula for the total symbol of the product of two
pseudodifferential operators, an $n$-dimensional differential form  $\Omega_n.$  
Note that we decided to 
write the constant on the left of \nnn{generalPolyakov} to simplify the typing.
This differential form is {\bf symmetric}, {\bf conformally 
invariant} and {\bf uniquely} 
determined, for every $f_0,f,h \in C^\infty(M),$ by the relation:
\begin{equation}\label{firsteq}
\Wres\bigl(f_0 [F,f][F,h]\bigr) =
\int_M f_0 \Omega_n(f,h).
\end{equation}

$\Omega_n(f,h)$ is given as the Wodzicki 1-density 
(which we denote $\wres$ following \cite{Polaris}) 
of the product of commutators $[F,f][F,h]$:
\begin{equation}\label{Omeganfh}
\begin{array}{rl}
&\Omega_n(f,h) = \wres([F,f][F,h])) \\
&\quad\quad=\biggl\{ 
   {\displaystyle \sum\frac{A_{\a',\a'',\b,\d}}{\a'!\a''!\b!\d!}} \, D^\b_x(f)D^{\a''+\d}_x(h)
\biggr\}\,d^{\,n}\, x,
\end{array}
\end{equation}
with the sum taken over $|\a'|+|\a''|+|\b|+|\d|+j+k=n,$ $|\b|\geq1,$ $|\d|\geq1,$ and
$$
A_{\a',\a'',\b,\d}=
{\displaystyle \int_{|\xi|=1}}\tr
\biggl\{
\del^{\a'+\a''+\b}_\xi\Bigl(\s_{-j}(F)\Bigr)
\del^{\d}_\xi\Bigl(D^{\a'}_x\Bigl(\s_{-k}(F)\Bigr)\Bigr)
\biggr\}\,d^{\,n-1}\,\xi
$$
where $\s_{-j}(f)=\s_{-j}(F)(x,\xi)$ is the component of order $-j$ in the total symbol of $F,$
$|\xi|=1$ means the Euclidean norm of the coordinate vector
$(\xi_1,\cdots,\xi_n)$ in $\bbR^n,$ and $d^{n-1}\xi$ is the normalized volume on
$\set{|\xi|=1}.$

In the 4-dimensional case, Connes \cite{Con1} showed that the Paneitz operator
$P_4,$
analogue of the scalar Laplacian in 4-dimensional conformal geometry, can be derived
from the Wodzicki residue by dropping some information.  That is to say, by setting 
$f_0 = 1$ and integrating by parts he has obtained
\begin{equation}
\label{Omega=fP}
\int_M \Omega_4(f,h) = \int_M f P_4(h) dv
\end{equation}
producing $P_4$ by the arbitrariness of $f$ and $h.$

The relation \nnn{Omega=fP} is of importance in the study of conformally invariant
differential operators generalizing the Yamabe operator.
There are (see \cite{GJMS}) invariant operators ({\em GJMS operators}) on scalar
densities with principal parts $\Delta^k,$ unless the dimension is even and $2k >n.$
The $n$-th order operator is called {\em critical} GJMS operator.  
When one calculates the Polyakov action for the quotient of functional determinants of a conformally
covariant operator $D$ at conformally related metrics, the operator $P_n$ shows, 
see for example \cite{Bra2} and \cite{BraO}.
That is why the study of such a differential 
form $\Omega_n,$ is of considerable importance in the case $n\geq 4.$  
We proposed an approach using automated symbolic computation to \nnn{Omeganfh}.

In this paper, we present partial results from \cite{Thesis}, see also \cite{paper2}.
In section~2 we introduce a general construction which associates to any 
pseudodifferential operator $S$ of order 0 acting on sections of a bundle $B$ on a compact manifold without 
boundary $M,$ a differential form of order $n$ acting on $C^\infty(M)\times C^\infty(M),$ $\Omega_{n,S}(f,h).$  
This $\Omega_{n,S}(f,h)$ is uniquely given by the relation 
$$
\Wres(f_0[S,f_1][S,f_2])=\int_M f_0 \Omega_{n,S}(f_1,f_2)\dvol
$$
for every $f_i \in C^\infty(M),$ with $\Wres$ the Wodsicki residue.  
The first result is given in 
Lemma~\nnn{lemmas_nFF} where we give an explicit expression for $\Omega_{n,S}$ in terms of the 
total symbol of $S.$
In the particular case of a compact, 
conformal, oriented, even-dimensional manifold, with $S=F$ the operator given by \nnn{defofF},
the differential form $\Omega_{n,F}$ is furthermore, symmetric on $f$ and $h$ and conformal invariant 
(Theorem~\nnn{theoremonOmegan}).  The rest of the paper focuses on computing $\Omega_{n,F}=\Omega_n$ 
in this case, in paticular for $n=6$.
In section~3, we give an explicit expression for $\Omega_n(f_1,f_2)$ in the flat case.  This 
expression (see Proposition~\nnn{propoomeganflat}) is given in terms of the Taylor expansion of the 
function $\tr(\s_L^F(\xi)\s_L^F(\eta))$ (see \nnn{trasFsF}).  
In section~4, we present $\Omega_6(f_1,f_2),$ given by Theorem~\nnn{theoremonOmegan}, in the 
6-dimensional conformally flat case.  
The last result is presented as Theorem~\nnn{lasttheorem} in section~6, where we
compute a whole family $\Omega_6(f,h)$ of differential forms of order 6 associated to 
a conformal, oriented, compact manifold without boundary.  Each of these forms will be symmetric 
on $f,$ and $h,$ conformally invariant and such that $\int_M f_0 \Omega_6(f_1,f_2)$ defines a 
Hochschild 2-cocycle over the algebra $C^\infty(M).$  To compute the unique form $\Omega_6$ 
satisfying the relation 
$\Wres(f_0[F,f_1][F,f_2])= \int_M f_0 \Omega_6(f_1,f_2)$ 
more information is needed in the 6-dimensional case.

My special thanks to T.~Branson for his constant support and guidance.  I am also in debt to the referee of 
this paper for valuable suggestions.

%%%###########%%% NOTATION SECTION %%%###########%%%

\medskip
\noindent{\bf{Notations and conventions}}

A conformal manifold is an equivalence class of Riemannian manifolds 
where two metrics $g$ and $\hat g$
are said to be equivalent if one is a positive scalar multiple of the other, for this work, it is
convenient to write $\hat g = e^{2\eta}g$ for some $\eta \in C^\infty(M).$

The Laplace-Beltrami operator on $k$-forms is
defined as $\Delta = d\d+\d d$ where we assume the sign convention $\Delta=-\frac{\del}{\del x} -
\frac{\del}{\del y}$ on $\bbR^2.$
The contraction of a $k$-form $\eta$ with a vector field $X$ is defined by
$$
\iota_\eta(X)(X_1,\cdots,X_{k-1}) := \eta(X,X_1,\cdots,X_{k-1}).
$$
The contraction of $\eta$ with a 1-form $\xi_X$ which is determined by the vector
field $X$ is given by $\iota_\eta(\xi_X) = \iota_\eta(X).$  The exterior multiplication
by a $k$-form $\eta$ will be denoted by $\eps_\eta : \xi \mapsto \eta\wedge\xi.$

In this work, the Riemann curvature tensor will be represented with the letter
$R,$ the Ricci tensor will be represented by $\Rc_{ij} = R^{k}{}_{ikj},$ and
the scalar curvature by $\Sc = \Rc^{i}{}_{i}.$ 
The conformal change equation for the Ricci tensor:
\begin{equation}
\label{invariantize}
\llcovd{\eta}{i}{j} =
- V_{ij}
- \lcovd{\eta}{i}\,\lcovd{\eta}{j}
+ \frac{1}{2}\lcovd{\eta}{k}\,\ucovd{\eta}{k}\,g_{ij},
\end{equation}
allows to replace, in the case $g=e^{2\eta} g_{\flat},$
the second derivatives on $\eta$ with terms with the Ricci tensor.
$V$ represents a normalized translation of the Ricci tensor,
useful in conformal geometry, given in terms of the normalized scalar curvature $J$ by
$$
V = \frac{\Rc - Jg}{n - 2}
\quad \hbox{with}\quad
J= \frac{\Sc}{2(n - 1)}.
$$
In \nnn{invariantize}, the indices after the semicolon represents covariant derivatives, 
$\llcovd{\eta}{i}{j}=\Na_j \Na_i \eta.$
In terms of $V,$ the relation between the Weyl tensor and the
Riemann tensor is given by
$$
W^i{}_{jkl} =
R^i{}_{jkl} + V_{jk}\d^i{}_l - V_{jl}\d^i{}_k + V^i{}_l g_{jk}
- V^i{}_k g_{jl}
$$
where $\d$ represents the Kronecker's delta tensor. 
  
If needed, we will ``raise'' and ``lower''
indices without explicit mention following 
% \cite{PenRin}, 
for example,
$g_{mi}R^i{}_{jkl}=R_{mjkl}.$
  
When working with the total symbol of a pseudodifferential operator 
$P,$ we will denote its leading symbol by $\s_L^P,$ or $\s_L(P)$ in case $P$ has a long expression.  
If the operator $P$ is of order $k$ then its total symbol (in some given local coordinates) 
will be represented as
$$
\s(P) = \s_k^P + \s_{k-1}^P + \s_{k-2}^P+\cdots,
$$
where $\s_L^P=\s_k^P.$  It is important to note that the different $\s_j^P$ for $j<k$ are defined 
only in local charts and are not diffeomorphism invariant \cite{Ku}.  However, Wodzicki \cite{Wo} has shown 
that the term $\s_{-n}^P$ enjoys a very special significance.   For a pseudodifferential operator $P,$ acting
on sections of a bundle $B$ over a manifold $M,$ there is a $1$-density on $M$ expressed in local coordinates
by
\begin{equation}
\label{wresP}
\wres(P) = \int_{|\xi|=1} \left\{\tr(\s_{-n}^P(x,\xi))\,d^{n-1}\xi\right\} \,d^n x.
\end{equation}
This \emph{Wodzicki residue density} is independent of the local representation.  Here we are using the same 
notations as in \cite{Polaris}, where an elementary proof of this matter can be found.
The Wodzicki residue, 
$\Wres(P),$ is then computed \cite{Con1} by choosing any local coordinates $x^j$ on $M$
and any local basis of sections $s_k$ for $B.$  $P$ is represented in
terms of the chosen basis $s_k$ as a matrix $P^l_k$ of scalar
psdeudodifferential operators: $P(f^k \a_k) =(P^i_k f^k)\a_i.$
The residue $\Wres(P)$ is given by
$$
\Wres(P^k_k)=
\int_M \left\{\int_{|\xi|=1} \tr(\s_{-n}(x,\xi))\,d^{n-1}\xi\right\}\,d^n x
$$
where $\s_{-n}(x,\xi)$ is the component of order $-n$ in the total symbol of $P,$
$|\xi|=1$ means the Euclidean norm of the coordinate vector
$(\xi_1,\cdots,\xi_n)$ in $\bbR^n,$ and $d^{n-1}\xi$ is the normalized volume on
$\set{|\xi|=1}.$
$\Wres(P)$ is independent of the choice
of the local coordinates on $M$, the local basis $(s_k)$ of $B$, and defines a trace (see \cite{Wo}).

To study the conformal invariance,
there is no need to study the whole conformal deformation. It is
enough to study the conformal deformation up to order one in $\eta$
as follows. If we set a metric $g$ in $M$ and consider another metric $\hat g$
conformally related to $g$ by the relation $\hat g = e^{2z\eta}g$ where
$\eta \in C^\infty(M)$ and $z$ a constant, then the conformal
variation of each expression is a polynomial in~$z$
whose coefficients are expressions in the metric and the conformal factor
$\eta$ (actually, this is an abuse of the language since the
conformal factor is $e^{2z\eta}$). In this way, the conformal
deformation up to order one in $\eta$ is given by
$\frac{d}{dz}\bigr|_{z=0}.$  As appointed in \cite{Gra}, if a natural tensor or a
differential operator is invariant up to order one in $\eta,$
i.e. if its conformal deformation up to order one is equal to zero, then by integration it
follows that it is fully invariant, for details see \cite{Bra4}.

%%%###########%%% END of NOTATIONS %%%###########%%%

%%%%%%%%%%%%%%%%%%%%%%%%section: Existence of $\Omega$  %%%%%%%%%%

\section{Existence of $\Omega_n$}

In this section, we associate to any pseudodifferential operator $S$ of order 0 acting on sections of a bundle $B$
on an arbitrary compact manifold without boundary $M$, a bilinear differential form $\Omega_{n,S}$ acting on 
$C^\infty(M)\times C^\infty(M)$.   Most of the properties of $\Omega_{n,S}$ are related with the properties of 
the Wodzicki residue. 
The total symbol up to order $-n$, of the pseudodifferential operator of
order $-2$ given by the product $P= f_0[S,f_1][S,f_2]$ with each $f_i \in C^\infty(M),$
is represented as a sum of $r\times r$ matrices of the form
$\s_{-2}^P+\s_{-3}^P+\cdots+\s_{-n}^P,$ with $r$ the rank of $B.$
We aim to study
$$
\Wres(P) = \int_M \left\{\int_{|\xi|=1}\tr(\s^P_{-n}(x,\xi))\,d^{n-1}\xi\right\}\,d^n x.
$$
  
In general, the total symbol of the product of two pseudodifferential
operators $P_1$ and $P_2$ is given by 
\begin{equation}
\label{sympro}
\sigma(P_1 P_2) = \sum \frac{1}{\a!}
\del^\a_\xi(\sigma^{P_1}) \,D^\a_x(\sigma^{P_2})
\end{equation}
where $\a=(\a_1,\dots,\a_n)$ is a multi-index,
$\a! = \a_1! \a_2! \cdots \a_n!$ and $D^\a_x = (-i)^{|\a|} \del^\a_x.$
Using this formula it is possible to deduce an expression for
$\s_{-n}([S,f_1][S,f_2])$ finding first $\s([S,f]).$

\begin{lemma}  
$[S,f]$ is a pseudodifferential operator of order $-1$ with total symbol
$\s([S,f]) = \sum_{k\geq 1} \s_{-k}([S,f])$ where
$$
\s_{-k}([S,f]) =\sum_{|\b|=1}^k \frac{1}{\b!} D^\b_x(f)
     \del_\xi^\b (\s^S_{-(k-|\b|)}).
$$
\end{lemma}

\begin{lemma}
\label{lemmas_nFF}
With the sum taken over $|\a'|+|\a''|+|\b|+|\d|+j+k = n,$ $|\b|\geq 1,$
and $|\d|\geq 1,$
\begin{align}
&\s_{-n}([S,f_1][S,f_2])
\nonumber \\
&=\sum \frac{1}{\a'!\a''!\b!\d!}
   D^\b_x(f_1) D_x^{\a''+\d}(f_2) \del_\xi^{\a'+\a''+\b}(\s^S_{-j})
\del_\xi^\d(D_x^{\a'}(\s^S_{-k})).
\label{sigma-nFF}
\end{align}
As a consequence
\begin{align}
&\wres([S,f_1][S,f_2]) 
= \biggl\{\int_{|\xi|=1} \tr\biggl\{ \sum \frac{1}{\a'!\a''!\b!\d!}
   D^\b_x(f_1) D^{\a''+\d}_x(f_2)
\times \nonumber \\
&\qquad \times
   \del^{\a'+\a''+\b}_\xi(\s^S_{-j}) \del^{\d}_\xi(D^{\a'}_x(\s^S_{-k}))
\biggr\}\,d^{n-1}\xi \biggr\}\,d^n x.
\end{align}
\end{lemma}

\begin{definition}
For every 
$f_1$ and $f_2$ in $C^\infty(M)$ we define
\begin{equation}
\label{defOmegan}
\Omega_{n,S}(f_1,f_2) := \wres([S,f_1][S,f_2]).
\end{equation}
\end{definition}

\begin{theorem}
For any pseudodifferential operator $S$ of order 0 acting on sections of a bundle $B$ on a manifold $M$, 
there is a unique $n$-diffe\-ren\-tial form $\Omega_{n,S}$ such that
$$
\Wres(f_0[S,f_1][S,f_2]) = \int_M f_0 \Omega_{n,S}(f_1,f_2)
$$
for all $f_i \in C^\infty(M).$
\end{theorem}

We restrict ourselves to an even dimensional, compact, oriented, conformal manifold without boundary $M$, and
$(B,S)$ given by the canonical Fredholm module $(\cH,F)$ associated to $M$ by A. Connes \cite{Con1}. 

In this particular case,
$F = (d\d - \d d)(d\d + \d d)^{-1}$ is the pseudodifferential operator of order 0
acting on the component $\Im(d+\d)$ of $\cH = L^2(M,\La^{l}_\bbC T^*M).$  
We can conclude that the differential form $\Omega_{n,F}$ is symmetric and conformally invariant.

\begin{theorem}
\label{theoremonOmegan}
In the particular case in which $M$ is a even dimensional compact conformal manifold without boundary and
$(\cH,F)$ is the Fredholm module associated to $M$ by A. Connes \cite{Con1},
there is a unique, symmetric, and conformally invariant $n$-diffe\-ren\-tial form $\Omega_n=\Omega_{n,F}$ 
such that
$$
\Wres(f_0[F,f_1][F,f_2]) = \int_M f_0 \Omega_n(f_1,f_2)
$$
for all $f_i \in C^\infty(M).$
\end{theorem}

All the proofs for the lemmas and theorems in this paper, as well as the detail for 
the computations in here presented can be read in \cite{Thesis}. 

%%%%%%%%%%%%%%%%%%%%section: The flat case %%%%%%%%%%%%%%%%%%%%%%%%%%%%%%%%%

\section{$\Omega_n$ in the flat case}

\begin{proposition}
\label{leadingsymbolofF}
The leading symbol of $F$ is given by
$$
\s_L^F(x,\xi) = \s_L^F(\xi) =
|\xi|^{-2}(\eps_\xi \iota_\xi - \iota_\xi \eps_\xi)
$$
for all $(x,\xi) \in T^* M,$ $\xi \not= 0.$ In the particular case of
a flat metric, we also have $\s_{-k}^F = 0$ for all $k \geq 1.$
\end{proposition}

Using this result, and the explicit expression for $\Omega_n$ (Definition~\ref{defOmegan}), it is possible to give a
formula for $\Omega_n$ in the flat case using the Taylor expansion of the function
\begin{equation}
\label{trasFsF}
\tr(\s_L^F(\xi) \s_L^F(\eta))=
a_{n,m} \frac{\<{\xi},{\eta}>^2}{|\xi|^2|\eta|^2} + b_{n,m}.
\end{equation}
Here  $\xi,\eta \in T^*_x M \smallsetminus \{0\},$ $\s_L^F(\xi) \s_L^F(\eta)$ is acting on $n/2$-forms,
and the values of the constants are given by
$$
\Biggl(\begin{array}{c}n \\ m\end{array}\Biggr) -  a_{n,m} =
b_{n,m} = \Biggl(\begin{array}{c}n-2 \\ m-2\end{array}\Biggr) +
          \Biggl(\begin{array}{c}n-2 \\ m\end{array}\Biggr) -
          2\Biggl(\begin{array}{c}n-2 \\ m-1\end{array}\Biggr).
$$

Because $\s_{-k}^F(x,\xi) = 0$
for all $k > 0$ in the flat case, \nnn{defOmegan} reduces to
\begin{align*}
&\Omega_{n\,\flat}(f,h) = \\
% &= \Bigl(\int_{|\xi|=1} \tr(\s_{-n}(x,\xi))\, d^{\,n-1}\xi \Bigr)\,\,d^{\,n}\,x \\
& \biggl\{ \int_{|\xi|=1} \tr\biggl( \sum \frac{1}{\a!\b!\d!}
(D^\b_x f)(D^{\a+\d}_x h)\del^{\a+\b}_\xi(\s_L^F)\del^\d_\xi(\s_L^F) \biggr)
 \,d^{n-1}\xi \biggr\} \,d^n x,
\end{align*}
with the sum taken over $|\a|+|\b|+|\d| = n,$ $1\leq |\b|,$ $1\leq |\d|.$

We denote by $T'_n\psi(\xi,\eta,u,v)$ the term of order $n$ in the
Taylor expansion of $\psi(\xi,\eta)=\tr(\s_L^F(\xi)\s_L^F(\eta))$ minus the terms with only powers
of $u$ or only powers of $v.$ That is to say,
$$
T'_n \psi(\xi,\eta,u,v) =
\sum_{|\b|+|\d|=n,|\b|\geq 1,|\d|\geq 1}
      \frac{u^\b}{\b!} \frac{v^\d}{\d!}
     \tr\bigl(\del^\b_\xi(\s_L^F(\xi)) \del^\d_\eta(\s_L^F(\eta))\bigr).
$$

\begin{proposition}
\label{propoomeganflat}
$$
\Omega_{n\,\flat}(f,h)
= \Bigl(\sum A_{a,b} (D_x^a f)(D_x^b h) \Bigr) \,d^n x,
$$
where
$$      
\sum A_{a,b} u^a v^b =
\int_{|\xi|=1} \bigl(T'_n\psi(\xi,\xi,u+v,v) -
                     T'_n\psi(\xi,\xi,v,v) \bigr) \,d^{n-1}\xi.
$$
\end{proposition}

%%%%%%%%%%%%%%%%%%%%%%%%%%%%%%%%%%%%%%%%%%%%%%%%%%%%%%%

%%%%%%%%%%%%%%%%%%%%section The 6-dimensional case  %%%%%%%%%%%%%%%%

\section{The 6-dimensional conformally flat case}

In the 6-dimensional case, symmetry and conformal invariance
are not enough to fully described $\Omega_6$ as we will find terms like
$$
\Bigl\{
A\,\lcovd{f}{i}\ucovd{h}{i}W_{jklm}W^{jklm}
+B\,\lcovd{f}{i}\lcovd{h}{j}W^i{}_{klm}W^{jklm}
\Bigr\}\,d^{\,6}\,x.
$$
which are symmetric on $f$ and $h,$ and conformally invariant.  

Other important property that 
we will exploit, is the fact that, by definition, $\Wres(f_0[F,f_1][F,f_2])$ is a Hochschild 2-cocycle 
over the algebra $C^\infty(M).$  If we define $\tau(f_0,f,h)$ as
$$
\tau(f_0,f,h)\colon
={\displaystyle \int_M}f_0\Bigl\{
A\,\lcovd{f}{i}\ucovd{h}{i}W_{jklm}W^{jklm}
+B\,\lcovd{f}{i}\lcovd{h}{j}W^i{}_{klm}W^{jklm}
\Bigr\}\,d^6\,x,
$$
then we obtain a Hochschild 2-cocycle for any value of $A$ and $B,$ that is to say \cite{Polaris}
\begin{align*}
0 &= (b\tau)(f_0,f_1,f_2,f_3) \\
  &= \tau(f_0f_1,f_2,f_3) - \tau(f_0,f_1f_2,f_3) +
     \tau(f_0,f_1,f_2f_3) - \tau(f_3f_0,f_1,f_2).
\end{align*}
In the 4-dimensional case, this property 
was used merely to make sure the constants found had the right values.  
In the 6-dimensional case, as we 
shall see, this property will play a more important role in the non-conformally flat case, even so, 
the fully description of $\Omega_6$ escapes these properties, requiring some more information to be used.

In the 6-dimensional flat case, using proposition~\ref{propoomeganflat} to compute $\Omega_6$ we have found
\begin{align}
&\Omega_{6\,\flat}(f,h)=Q_6(df,dh) \nonumber\\
&=\Bigl\{12(\lcovd{f}{i}\ululucovd{h}{i}{j}{j}{k}{k} +
                  \llulucovd{f}{i}{j}{j}{k}{k}\ucovd{h}{i})
             + 24\,(\llcovd{f}{i}{j}\uulucovd{h}{i}{j}{k}{k} +
                 \lllucovd{f}{i}{j}{k}{k}\uucovd{h}{i}{j}) \nonumber\\
&\quad      + 6\,(\lucovd{f}{i}{i}\lulucovd{h}{j}{j}{k}{k} +
                 \lulucovd{f}{i}{i}{j}{j}\lucovd{h}{k}{k})
             + 24\,\llucovd{f}{i}{j}{j}\ulucovd{h}{i}{k}{k}
             + 16\,\lllcovd{f}{i}{j}{k}\uuucovd{h}{i}{j}{k}\Bigr\}\, d^{\,6}\,x \nonumber\\
&=\Bigl\{12\,\Delta^2(\<df,dh>)-6 \,\Delta(\Delta f \Delta h) 
             - 12\,\<\Na(\Delta f),\Na(\Delta h)>\nonumber\\
&\quad  +24\, \Delta(\<\Na df,\Na dh>)+16\, \<\Na^2 df,\Na^2 dh> \Bigr\} d^{\,6}\, x,
\label{omegaflat}
\end{align}
where each summand in the last expression is explicitly symmetric on $f$ and $h.$

Studying $\Omega_n$ in the flat metric gives information
about the conformally flat case, in particular, by growing the flat expression
for $\Omega_n$ we can find its expression in the conformally flat case.
To do that, we consider a metric $\hat g$ conformally related to the flat metric $g$ by the relation
$\hat g=e^{2\eta}g$ with $\eta$ an smooth function on $M.$
Each time we express a component of \nnn{omegaflat} in the conformally related metric, terms containing
derivatives on $\eta$ will show.  Using the conformal change 
equation for the Ricci tensor \nnn{invariantize} until 
we reduce all the higher derivatives on $\eta$ to derivatives of order one,   
the Ricci tensor, and the scalar curvature related to~$g,$ we obtain the expression for $\Omega_n$ 
in the conformally flat case.
The computations could be a little tedious, so we did 
them using \texttt{Ricci.m} \cite{Lee} to obtain: 
\begin{equation}
\label{omega6confflat}
\begin{array}{rl}
&\Omega_{6\,\confflat}(f,h)= \Omega_{6\,\flat}(f,h) \\
& + \Bigl\{\Bigl(
         -72(\llucovd{f}{i}{j}{j}\ucovd{h}{i} + \lcovd{f}{i}\ulucovd{h}{i}{j}{j})
 -24\lucovd{f}{i}{i}\lucovd{h}{j}{j} - 96\llcovd{f}{i}{j} \uucovd{h}{i}{j}\Bigr)J \\ 
& + 96\lcovd{f}{i}\ucovd{h}{i} J^2 \\
& + 24(\lucovd{f}{i}{i}\lcovd{h}{j} \ucovd{J}{j}
       + \lcovd{f}{i}\lucovd{h}{j}{j} \ucovd{J}{i}) 
       - 24(\llcovd{f}{i}{j}\ucovd{h}{i}\ucovd{J}{j}
          + \lcovd{f}{i}\ulcovd{h}{i}{j}\ucovd{J}{j})\\
& - 24\lcovd{f}{i}\ucovd{h}{i}\lucovd{J}{j}{j}
       + 64 \lcovd{f}{i}\lcovd{h}{j} J V^{ij}
 - 32(\llucovd{f}{i}{j}{j}\lcovd{h}{k}V^{ik}
           + \lcovd{f}{i} \llucovd{h}{j}{k}{k} V^{ij})\\
& +64(\lllcovd{f}{i}{j}{k}\ucovd{h}{i}V^{jk}
                    + \lcovd{f}{i} \ullcovd{h}{i}{j}{k} V^{jk})
       +96(\llcovd{f}{i}{j}\lucovd{h}{k}{k} V^{ij} 
                    + \lucovd{f}{i}{i} \llcovd{h}{j}{k} V^{jk})\\
& -192\llcovd{f}{i}{j} \ulcovd{h}{j}{k} V^{jk}
       -64 \lcovd{f}{i}\ucovd{h}{i} V_{jk}V^{jk} 
       +128 \lcovd{f}{i}\lcovd{h}{j}V^i{}_{k}V^{jk}\Bigr\}d^{\,6}\, x.
\end{array}
\end{equation}

An interesting introduction to automated symbolic computations 
can be found in \cite{Bra1}.

\begin{theorem}
In the 6-dimensional conformally flat case,
the expression for $\Omega_6$ given by Theorem~\nnn{theoremonOmegan} 
as a sum of explicitly symmetric components on $f$ and $h,$ is given by 
\begin{equation}
\label{omega6confflatsym}
\begin{array}{rl}
&\Omega_{6 \confflat}(f,h)
=\Bigl\{12\, \Delta^2(\<df,dh>)-6\Delta(\Delta f \Delta h) - 12\<\Na \Delta f,\Na \Delta h> \\
& + 24\, \Delta(\<\Na df,\Na dh>)+16\, \<\Na^2 df,\Na^2 dh> +72\, \Delta\<df,dh>J \\
& - 24\Delta(f)\Delta(h)J + 48\,\<\Na df,\Na dh>J + 96\,\<df,dh>J^2 + \\
& + 24\,\<df,dh> \Delta(J) -24\<\Delta(f)dh + \Delta(h)df,dJ> \\
& - 24\, \<{d(\<df,dh>)},dJ> - 96\<\Delta(h){\Na df} + \Delta(f){\Na dh},V> \\
& + 32\, \<{\Na(\Delta(f)\otimes dh)} + {\Na(\Delta(h)\otimes df)},V> \\
& + 64\,\<{\Na^2(\<df,dh>)},V> -64\, \<df,dh>\<V,V> \\
& - 128\, \tr((df \otimes dh) V^2) + 64 \, \tr((\Hess f)(\Hess h)V)   \Bigr\} d^{\,6}\, 
x.
\end{array}
\end{equation}
In the last two terms, both factors are considered
as $(1,1)$ tensors (one contravariant and one covariant).
\end{theorem}

Actually, the difference in between the two expressions \nnn{omega6confflat}
and \nnn{omega6confflatsym} is given by the term
$$
\nnn{omega6confflatsym}-\nnn{omega6confflat}
=
96\,\llcovd{f}{i}{j}\llcovd{h}{k}{l} W^{iljk}
-32(\llcovd{f}{i}{j}\lcovd{h}{k}\ucovd{W^{ijk}{}_l}{l}+
    \lcovd{f}{i}\llcovd{h}{j}{k}\ucovd{W^{ijk}{}_l}{l})
$$
which vanishes in the conformally flat case. 
 
Leaving for an instant the conformally flat case, in the general conformally curved case,
the conformal variation of $\Omega_6(f,h),$ up to order one in $\eta$ is given by
\begin{equation}\label{conformalvaritionomega6cf}
\begin{array}{rl}
& -32
  \Bigl\{\lcovd{\eta}{i} \lcovd{f}{j} \lcovd{h}{k} \ucovd{W^{ijk}{}_{l}}{l}
    +\lcovd{\eta}{i} \lcovd{f}{j} \lcovd{h}{k} \ucovd{W^{ikj}{}_{l}}{l} \\
&\quad\quad\quad +\lcovd{\eta}{i}\llcovd{f}{j}{k} \lcovd{h}{l}W^{ijkl}
        - \lcovd{\eta}{i}\lcovd{f}{j}\llcovd{h}{k}{l} W^{ikjl}\Bigr\}\,d^{\,6}\,x,
\end{array}
\end{equation}
which vanishes in the conformally flat case, meaning that our 
expression is conformally invariant inside the conformally flat class of metrics on $M.$ In 
the general conformally curved case, this variation will be useful in finding the extra 
terms we are missing, that is to say, those terms that vanish in the conformally flat case.

If we define using \nnn{omega6confflat} the trilinear form on $C^\infty(M)$ 
$$
\tau(f_0,f_1,f_2):= \int_M f_0\Omega_{6,\confflat}(f_1,f_2)
$$
then
\begin{align}
&(b \tau)(f_0,f_1,f_2,f_3) \nonumber\\
&=\int_M f_0 
\Bigl(-96\,(
\lcovd{f_1}{j}\lcovd{f_2}{i}\lcovd{f_3}{k}\ucovd{W^{ijk}{}_l}{l}
            +\lcovd{f_1}{j}\lcovd{f_2}{i}\lcovd{f_3}{k}\ucovd{W^{ikj}{}_l}{l})
\nonumber\\
&\quad\quad\quad\quad+128\,(\llcovd{f_1}{j}{k}\lcovd{f_2}{i}\lcovd{f_3}{l}W^{ijkl}+
       \lcovd{f_1}{j}\lcovd{f_2}{i}\llcovd{f_3}{k}{l}W^{ikjl})\Bigr)\,d^6x
\label{omconflaH}
\end{align}
which vanishes in the conformally flat case meaning that $\tau$ is a Hochschild 2-cocycle, 
in the conformally flat case.
 
%%%%%%%%%%%%%%%%%%%%%%%%%%%%%%%%%%%%%%%%%%%%%%%%%%%%%%%%%%%%%%%%%%%%%

\section{A filtration by degree}
To simplify the notation, and because of the factor $d^n x$ in
the definition of $\Omega_n$ we will write
$\Omega_n(f,h) = \omega_n(f,h)\,d^n x$ where,
\begin{align*}
\omega_n(f,h)
&= \sum \frac{1}{\a'!\a''!\b!\d!} D^\b_x(f) D^{\a''+\d}_x(h) \times\\
&\qquad \times \int_{|\xi|=1}
\tr\bigl(\del^{\a'+\a''+\b}_\xi(\s^F_{-j})
         \del^{\d}_\xi(D^{\a'}_x(\s^F_{-k}))\bigr) \,d^{n-1}\xi.
\end{align*}
The expression $\omega_n(f,h)$ is a sum of homogeneous polynomials 
in the ingredients
$\Na^\a df,$ $\Na^\b dh,$ and $\Na^\ga R$ for multiindices $\a,\b,$ and
$\ga,$ in the following sense, each
monomial must satisfies the {\bf homogeneity condition}  given by the
rule (see \cite{Bra4}):
$$
\mbox{twice~the~appearances~of~$R$~}
+ \mbox{~number~of~covariant~derivatives} = n
$$
where for covariant derivatives we count all of the derivatives on
$R,$ $f,$ and $h,$ and any occurrence of $W,$ $Rc,$ $V,$ $Sc,$ or $J$
is counted as an occurrence of~$R.$
By closing under addition, we denote by $\cP_n$ the space of these
polynomials.

This same idea is used in \cite{BraGiOr} to study leading terms in the
heat invariants produced by the Laplacian of de~Rham and other
complexes. We borrow from there the idea of \emph{filtration by degree}.
For
a homogeneous polynomial $Q$ in $\cP_n,$ we denote by $k_R$ its degree
in~$R$ and by $k_\Na$ its degree in~$\Na.$  In this way,
$2k_R + k_\Na = n.$
Because $|\b|\geq 1,$ and $|\d|\geq 1$ in Lemma~\ref{lemmas_nFF}, we
have $k_\Na \geq 2$ and hence $2k_R \leq n-2.$

We say that $Q$ is in $\cP_{n,l}$ if $Q$ can be written as a sum of
monomials with $k_R\geq l,$ or equivalently, $k_\Na \leq n - 2l.$ We
have
$$
\cP_n = \cP_{n,0} \supseteq \cP_{n,1} \supseteq \cP_{n,2}
\supseteq\cdots\supseteq \cP_{n,\tfrac{n-2}{2}},
$$
and $\cP_{n,l} = 0$ for $l > (n-2)/2.$
There is an important observation to make. An expression which a
priori appears to be in, say $\cP_{6,1},$ may actually be in a
subspace of it, like $\cP_{6,2}.$  For example,
$$
\underbrace{\lcovd{f}{i} \lllcovd{h}{j}{k}{l} W^{ijlk}}_{\in\,\cP_{6,1}} =
\underbrace{\lcovd{f}{i} \lcovd{h}{j} V_{kl} W^{ikjl}
+\lcovd{f}{i}\lcovd{h}{j} W^i{}_{klm} W^j{}_{klm}}_
{\in \, \cP_{6,2}},
$$
by reordering covariant derivatives and making use of the symmetries of
the
Weyl tensor.
Because of this filtration, we use a fix convention on how the indices
should
be placed when representing each expression in its index notation.  
For example, $\lllcovd{f}{i}{j}{k} \uuucovd{h}{i}{j}{k}$ will be
preferred over $\lllcovd{f}{i}{k}{j} \uuucovd{h}{i}{j}{k}.$
Also$\llulucovd{f}{i}{j}{j}{k}{k}$ will be preferred
over $\lullucovd{f}{i}{i}{j}{k}{k}$ or any other variation.
Once we have defined the filtration on $\cP_n,$ and accepted our
notational convention, we can state the following proposition.

\begin{proposition}
\label{omeganfiltration}
There exists a universal bilinear form $Q_n(df,dh)$ in \newline
$\cP_{n,0} \smallsetminus \cP_{n,1}$ and a form $Q_{R,n}(df,dh)$ in
$\cP_{n,1}$ such that 
$$
\Omega_n(f,h) = Q_n(df,dh) + Q_{R,n}(df,dh).
$$
In the particular case of the flat metric $\Omega_n(f,h) = Q_n(df,dh)$
since the curvature vanishes.
\end{proposition}

In the particular case $n=4,$ $k_R$ can be 0 or~1, hence $\Omega_4$
can be written as
$$
\Omega_4(f,h) = \underbrace{Q_4(df,dh)}_{\in\,\cP_{4,0}} +
\underbrace{Q_{R,4}(df,dh)}_{\in\,\cP_{4,1}}
$$
where $Q_{R,4}(df,dh)$ is a trilinear form on $R,$ $df,$ and $dh.$

In the 6-dimensional case,
$k_R \in \{0,1,2\}$ thus
\begin{align}
\Omega_6(f,h) &= Q_6(df,dh)
\nonumber \\
&\qquad + Q_{R,6}^{(1,0)}(df, dh)
        + Q_{R,6}^{(1,1)}(df, dh) + Q_{R,6}^{(1,2)}(df, dh)
\nonumber \\
&\qquad + Q_{R,6}^{(2,0)}(df, dh),
\label{omega6filtration}
\end{align}
where
\begin{itemize}  
\item[-]
$Q_{R,6}^{(1,0)}(df, dh) \in \cP_{6,1} \smallsetminus \cP_{6,2},$  
without covariant derivatives on $R,$
\item[-]
$Q_{R,6}^{(1,1)}(df, dh) \in \cP_{6,1} \smallsetminus \cP_{6,2},$
with a single covariant derivative on $R,$  
\item[-]
$Q_{R,6}^{(1,2)}(df, dh) \in \cP_{6,1} \smallsetminus \cP_{6,2},$ 
with two covariant derivatives on $R,$ and
\item[-]
$Q_{R,6}^{(2,0)}(df, dh) \in \cP_{6,2},$ without covariant
derivatives on $R.$
\end{itemize}

{}From the previous expressions, it is evident that there exists a
sub-filtration inside each $\cP_{n,l}$ for $l\geq 1.$ Such a
filtration is a lot more complicated to describe in higher dimension
because of the presence of terms like $\Na^a R \Na^c R\cdots.$

%%%%%%%%%%%%%%%%%%%%%%%%%%%%%%%%%%%%%%%%%%%%%%%%%%%%%%%%%%%%%%%%%%%%%

\section{The 6-dimensional non-conformally flat case}

We do not restrict ourselves anymore to the conformally flat case.  Now we are going to find those 
terms we need to add to $\Omega_{6,\confflat}$ in order to get the expression for $\Omega_6$ in the 
general conformally curved case.

The first set of terms to be added will complete the expression for 
$Q_{R,6}^{(1,0)}.$  In this case, there is just one possibility to be consider, that is
$$
\llcovd{f}{i}{j}\llcovd{h}{k}{l} W^{ikjl}.
$$  
Any other possibility is ruled out by the relation
$$
\lcovd{f}{i} \lllcovd{h}{j}{k}{l} W^{ijkl} = 
\lcovd{f}{i} \lcovd{h}{j} V_{kl} W^{ikjl} 
+\lcovd{f}{i}\lcovd{h}{j} W_i{}^{klm} W^{jlkm}
$$
which express the right hand side, an element of $\cP_{6,1},$ as the sum of two elements of 
$\cP_{6,2}.$

The second set of terms completes the expression for 
$Q_{R,6}^{(1,1)},$ it is given by
the following symmetric term on $f$ and $h$:
$$
\llcovd{f}{i}{j}\lcovd{h}{k} \ucovd{W^i{}_l{}^{jk}}{l} 
+ \llcovd{h}{i}{j}\lcovd{f}{k} \ucovd{W^i{}_l{}^{jk}}{l}.
$$

The only term to complete the expression for $Q_{R,6}^{(1,2)}$ is
$$
\lcovd{f}{i} \lcovd{h}{j} \uucovd{W^i{}_k{}^j{}_l}{k}{l},
$$
symmetric on $f$ and $h.$

For $Q_{R,6}^{(2,0)}$ we consider at this time, just one term
$$
\lcovd{f}{i}\lcovd{h}{j} V_{kl}W^{ikjl}.
$$
It happens that the other possible terms 
\begin{equation}\label{badterms}
\lcovd{f}{i}\ucovd{h}{i}W_{jklm}W^{jklm}\,d^6x\quad\hbox{and}\quad
\lcovd{f}{i}\lcovd{h}{j}W^{i}{}_{klm}W^{jklm}\,d^6x.
\end{equation}
are conformally invariant.  

Up to this point, what we must add to $\Omega_{6,\confflat}$ is a linear combination of the form
$$
\begin{array}{rl}
&\Bigl\{A\,\llcovd{f}{i}{j}\llcovd{h}{k}{l} W^{ikjl}+
B\,(\llcovd{f}{i}{j}\lcovd{h}{k} \ucovd{W^i{}_l{}^{jk}}{l} 
+ \llcovd{h}{i}{j}\lcovd{f}{k} \ucovd{W^i{}_l{}^{jk}}{l} )\\
&\quad + C\,\lcovd{f}{i} \lcovd{h}{j} \uucovd{W^i{}_k{}^j{}_l}{k}{l}
       + D\,\lcovd{f}{i}\lcovd{h}{j} V_{kl}W^{ikjl}\Bigr\}\,d^{\,6}\,x.
\end{array}
$$

Its conformal variation up to order one in $\eta$ is given by 
$$
\begin{array}{rl}
&\Bigl\{(B+2C)\,(\lcovd{\eta}{i}\lcovd{f}{j}\lcovd{h}{k}\ucovd{W^{ijk}{}_l}{l}
           +\lcovd{\eta}{i}\lcovd{f}{j}\lcovd{h}{k}\ucovd{W^{ikj}{}_l}{l})\\
&\quad+(3B-2A)\,(\lcovd{\eta}{i}\llcovd{f}{j}{k}\lcovd{h}{l} W^{ijkl} 
              + \lcovd{\eta}{i}\lcovd{f}{j}\llcovd{h}{k}{l} W^{ilkj})\\
&\quad+(D-3C)\,(\llcovd{\eta}{i}{j}\lcovd{f}{k}\lcovd{h}{l}W^{ikjl}) \Bigr\}\,d^6x.
\end{array}
$$

By comparing it with the conformal variation of $\Omega_{6,\confflat}(f,h)$ 
\nnn{conformalvaritionomega6cf}, we deduce the conditions $B+2C=-32=3B-2A,$ 
and $D-3C=0,$ which means, conformally invariance and symmetry is not enough to find the right 
values for all the constants.  So far, the term to be added to $\Omega_6(f,h)$ is given by
\begin{equation}
\label{termstoadd1}
\begin{array}{rl}
&\Bigl\{-(32+3C)\,\llcovd{f}{i}{j}\llcovd{h}{k}{l} W^{ikjl}\\
&-(32+2C)\,(\llcovd{f}{i}{j}\lcovd{h}{k} \ucovd{W^i{}_l{}^{jk}}{l} 
+ \llcovd{h}{i}{j}\lcovd{f}{k} \ucovd{W^i{}_l{}^{jk}}{l} \\
&\quad + C\,\lcovd{f}{i} \lcovd{h}{j} \uucovd{W^i{}_k{}^j{}_l}{k}{l}
       + 3C\,\lcovd{f}{i}\lcovd{h}{j} V_{kl}W^{ikjl}\\
&\quad + E\,\lcovd{f}{i}\ucovd{h}{i}W_{jklm}W^{jklm}
       + G\,\lcovd{f}{i}\lcovd{h}{j}W^{i}{}_{klm}W^{jklm}\Bigr\}\,d^{\,6}\,x
\end{array}
\end{equation}
where the last two terms come from \nnn{badterms}.

%%%%%%%%%%%%%%%%%  Subsection: Using the Hoschild 2-cocycle property  %%%%%%%%%%%%%

\medskip
\noindent{{\bf Using the Hochschild 2-cocycle property}}

\begin{proposition} 
If we define the trilinear form on $C^\infty(M)$
\begin{align*}
&\tau'(f_0,f_1,f_2):= \\
&\int_M f_0
\Bigl\{C\,\lcovd{f_1}{i} \lcovd{f_2}{j} \uucovd{W^i{}_k{}^j{}_l}{k}{l}
       + D\,\lcovd{f_1}{i}\lcovd{f_2}{j} V_{kl}W^{ikjl}\\
& + E\,\lcovd{f_1}{i}\ucovd{f_2}{i}W_{jklm}W^{jklm}
       + G\,\lcovd{f_1}{i}\lcovd{f_2}{j}W^{i}{}_{klm}W^{jklm}\Bigr\}\,d^{\,6}\,x
\end{align*}
then 
$$
(b\tau')(f_0,f_1,f_2,f_3)=0
$$ 
for any $f_i \in C^\infty(M)$ meaning that we obtain a Hochschild 2-cocycle on the 
algebra $C^\infty(M)$ for any value of the constants $C,$ $D,$ $E,$ and $G.$  
\end{proposition}

On the other hand, if we define
\begin{align*}
\tau''(f_0,f_1,f_2):= \int_M f_0
&\Bigl\{A\,\llcovd{f_1}{i}{j}\llcovd{f_2}{k}{l} W^{ikjl} \\
&\quad+B\,(\llcovd{f_1}{i}{j}\lcovd{f_2}{k} \ucovd{W^i{}_l{}^{jk}}{l}
+ \llcovd{f_2}{i}{j}\lcovd{f_1}{k} \ucovd{W^i{}_l{}^{jk}}{l} )\Bigr\}\,d^{\,6}\,x
\end{align*}
then
\begin{align}
&(b \tau'')(f_0,f_1,f_2,f_3) \nonumber\\
&=\int_M f_0
\Bigl( 3\,B\,(\lcovd{f_1}{j}\lcovd{f_2}{i}\lcovd{f_3}{k}\ucovd{W^{ijk}{}_l}{l}
            +\lcovd{f_1}{j}\lcovd{f_2}{i}\lcovd{f_3}{k}\ucovd{W^{ikj}{}_l}{l})
\nonumber\\
&\quad\quad\quad\quad
-2\,A\,\,(\llcovd{f_1}{j}{k}\lcovd{f_2}{i}\lcovd{f_3}{l}W^{ijkl}+
       \lcovd{f_1}{j}\lcovd{f_2}{i}\llcovd{f_3}{k}{l}W^{ikjl})
\Bigr)\,d^6x.
\label{termstoaddH}
\end{align}

To have that $\int_M f_0\Omega_6(f_1,f_2)$ is a Hochschild 2-cocycle on $C^\infty(M)$ we 
need $\nnn{omconflaH} +\nnn{termstoaddH}=0,$ for any $f_i \in C^\infty(M).$  
Thus $3B=96$ and $2A=128.$  Because $B+2C=-32=3B-2A$ we must have $C=-32$ and hence 
using \nnn{omega6confflat}, \nnn{badterms}, and \nnn{termstoadd1} we conclude
\begin{align}
\Omega_6(f,h)&=\Omega_{6\,\confflat}(f,h)\nonumber\\
&\quad+ \Bigl\{64\,\llcovd{f}{i}{j}\llcovd{h}{k}{l} W^{ikjl}
+ 32\,(\llcovd{f}{i}{j}\lcovd{h}{k} \ucovd{W^i{}_l^{jk}}{l} 
+ \llcovd{h}{i}{j}\lcovd{f}{k} \ucovd{W^i{}_l^{jk}}{l}) \nonumber\\
&\quad\quad -32\,\lcovd{f}{i} \lcovd{h}{j} \uucovd{W^i{}_k{}^j{}_l}{k}{l}
       -96\,\lcovd{f}{i}\lcovd{h}{j} V_{kl}W^{ikjl}\nonumber\\
&\quad\quad+E\,\lcovd{f}{i}\ucovd{h}{i}W_{jklm}W^{jklm}
+G\,\lcovd{f}{i}\lcovd{h}{j}W^{i}{}_{klm}W^{jklm}\Bigr\}\,d^{\,6}\,x,
\label{last}
\end{align}
where the last two terms are the needed ones to fully complete the expression for
$Q^{(2,0)}_{R,6}$ as in \nnn{badterms}.

\begin{theorem}
\label{lasttheorem}
The expression \nnn{last} gives a family of 6-dimensional differential forms associated to $M,$ 
each of these differential forms is symmetric on $f$ and $h,$ conformally invariant, and defines 
a Hochschild 2-cocycle on the algebra $C^\infty(M)$ by the relation 
$\tau(f_0,f_1,f_2)=\int_M f_0\Omega_6(f_1,f_2).$
\end{theorem}

In the 6-dimensional conformally curved case, more information is needed to find the unique one satisfying 
the relation 
$$
\Wres(f_0[F,f_1][F,f_2])=\int_M f_0\Omega_6(f_1,f_2).
$$

%%%%%%%%%%%%%%%%%%%%%%%%%%%%%%%%%%%%%%%%%%%%%%%%%%%%%%%%%%%%%%%%

%%%%%%%%%%%%%%%%%%%%%information about the author
\small
\vskip 1pc
{\obeylines
\noindent William \ J. \ Ugalde 
\noindent Department of Mathematics
\noindent 14 MacLean Hall
\noindent The University of Iowa
\noindent Iowa City, Iowa, 52242
\noindent E-mail: wugalde@math.uiowa.edu
}


\begin{thebibliography}{99}
\def\topsep{0pthere exists a}
\def\parsep{0pt plus 5pt minus 1pt}
\def\itemsep{-0.5ex} %seems to be a nice small skip between items
\small               %seems to be the best
%
\bibitem{Bra1}
T. P. Branson,
{\it Automated symbolic computations in spin geometry\/},
In ``Clifford Analysis and its Applications,'' F.\ Brackx, J.S.R.\
Chisholm, and V.\ Sou\v cek, eds.,
pp.\ 27--38.  NATO Science Series II, Vol.\ 25,
Kluwer Academic Publishers, Dordrecht, 2001.
%
\bibitem{Bra4} T. P. Branson,
{\it Differential operators canonically associated to a conformal structure\/},
Math. Scand., \textbf{57} (1985), 293--345.
%
\bibitem{Bra2} T. Branson, {\it An Anomaly associated with 4-dimensional quantum
gravity\/}, Commun. Math. Phys. 178 (1996) 301--309.
%
\bibitem{BraGiOr} T.P. Branson, P.B. Gilkey, and Bent \O rsted, {\it Leading terms in the
heat invariants for the Laplacians of the de Rham, signature, and spin complexes\/}, 
Math. Scand. 66 (1990), 307--319.
%
\bibitem{BraO} T. Branson and B. \O rsted, {\em Explicit functional determinants in
four dimensions}, Proc. Amer. Math. Soc. 113 (1991) 669--682.
%
\bibitem{Con1} A. Connes, {\it Quantized calculus and applications\/}, XIth International
Congress of Mathematical Physics (Paris, 1994), 15--36, Internat. Press, Cambridge, MA,
1995.
%
\bibitem{Koran}
A. Connes,
{\it Noncommutative Geometry\/},
Academic Press, London and San Diego, 1994.
%
\bibitem{GoPe} A. R. Gover, L. J. Peterson, {\it Conformally invariant powers of the
Laplacian, Q-curvature and tractor calculus\/}, to appear.
%
\bibitem{Polaris}
J. M. Gracia-Bond\'{\i}a, J. C. V\'arilly and H. Figueroa,
{\it Elements of Noncommutative Geometry\/},
Birkh\"auser Advanced Texts, Birkh\"auser, Boston, 2001.
%
\bibitem{Gra}
C. R. Graham,
{\it Conformally invariant powers of the Laplacian, II: Nonexistence\/},
J. London Math. Soc. (2) \textbf{46} (1992), 566--576.
%
\bibitem{GJMS} R. Graham, R. Jenne, L. Mason, and G. Sparling, {\em Conformally
invariant powers of the Laplacian, I: Existence}, J. London Math. Soc. (2) 46 (1992)
557--565.
%
\bibitem{Ku}
H. Kumano-go, {\it Pseudo-Differential Operators\/},
The MIT Press, Cambridge, Massachusetts, 1981.
%
\bibitem{Lee}
J. M. Lee, {\it A Mathematica package for doing tensor calculations in
differential geometry\/}, http://www.math.washington.edu/~lee/Ricci/
%
% \bibitem{PenRin}
% R. Penrose and W. Rindler,
% {\it Spinors and Space-time, vol.1\/},
% Cambridge University Press, Cambridge, 1984.
%
\bibitem{Thesis}
W. J. Ugalde, {\it On differential forms canonically associated to even dimensional 
conformal manifolds\/}, The University of Iowa, in progress.
%
\bibitem{paper2}
W. J. Ugalde, {\it Differential forms and the Wodzicki residue\/}, arXiv:math-DG/0211361.
%
\bibitem{Warner}
F. W. Warner, {\it Foundations of Differential Manifolds and Lie Groups\/},
Springer-Verlag, New York, 1983.
%
\bibitem{Wo}
M. Wodzicki, {\it Noncommutative residue, Part I.  Fundamentals, $K$-Theory, arithmetic
and geometry\/}, (Moscow 1984-86), pp. 320-399, Lecture Notes Math, 1289,
Springer, Berlin, 1987.

\end{thebibliography}
\end{document}